\documentclass[a4paper,10pt]{article}
\usepackage[utf8x]{inputenc}
\usepackage{amsthm,amsmath,amssymb,oldgerm}
\usepackage[margin=3cm]{geometry}
\theoremstyle{plain}
\newtheorem{thm}{Theorem}[section]
\newtheorem*{theorem*}{Main theorem}

\theoremstyle{definition}
 
\theoremstyle{definition}
\newtheorem{rem}[thm]{Remark}

\let\H\relax
\let\L\relax
\let\O\relax
\let\dh\relax
\DeclareMathOperator{\cx} {\mathcal{C}_{{\it{X}}}}
\DeclareMathOperator{\bk} {\mathcal{B}_{{\it{X}}}^{2{\it{k}}}}

\DeclareMathOperator{\bkx} {\mathcal{B}_{\Omega_{\it{X}}}^{\it{k}}}
\DeclareMathOperator{\hyp}{\mu_{\mathrm{hyp}}} 
 
\DeclareMathOperator{\H}{\mathbb{H}}
\DeclareMathOperator{\N}{\mathbb{N}}  
\DeclareMathOperator{\G}{\Gamma} 
\DeclareMathOperator{\L}{\mathcal{L}} 
\DeclareMathOperator{\C}{\mathbb{C}}
\DeclareMathOperator{\R}{\mathbb{R}} 
 
\DeclareMathOperator{\O}{\Omega_{{\it{X}}}} 
\DeclareMathOperator{\Ok}{\Omega_{{\it{X}}}^{\otimes{\it{k}}}} 
\DeclareMathOperator{\Stwo}{\mathcal{S}^{2}(\Gamma)} 
\DeclareMathOperator{\Sk}{\mathcal{S}^{2{\it{k}}}(\Gamma)}
\DeclareMathOperator{\jk}{{\it{j}}_{{\it{k}}}} 
\DeclareMathOperator{\rx}{{\it{r}}_{{\it{X}}}} 
\DeclareMathOperator{\dh}{{\it{d}}_{\mathbb{H}}}
\title{Off-diagonal estimates of the Bergman kernel on hyperbolic Riemann surfaces of finite volume}
\author{Anilatmaja Aryasomayajula and Priyanka Majumder}
\date{}
\begin{document}
\maketitle
\begin{abstract}
\vspace{0.2cm}\noindent
In this article, we derive off-diagonal estimates of the Bergman kernel associated to tensor-products of the cotangent line bundle defined over a hyperbolic Riemann surface of finite volume. 

\vspace{0.1cm}\noindent
\end{abstract}
\section{Introduction}
In this article, we derive off-diagonal quantitative estimates of the Bergman kernel associated to tensor-powers of the cotangent bundle defined on a hyperbolic Riemann surface of finite volume. 

\vspace{0.1cm}
Estimates of Bergman kernels associated to high tensor-powers of holomorphic line bundles defined on complex manifolds has been an object of study since a long time. Tian, Zelditch, Demailly, Marinsecu, Ma et al. have done seminal work in this field. 

\vspace{0.1cm}
In this article, we derive estimates of the Bergman kernel associated to tensor-powers of the cotangent bundle defined over a hyperbolic Riemann surface of finite volume, away from the diagonal, both in the compact and in the noncompact setting. Our estimates depend only on the injectivity radius of the hyperbolic Riemann surface, and the tensor-power of the cotangent bundle. 

\vspace{0.1cm}
{\bf{Results from Literature.}} We now briefly discuss the history behind the problem, before we state our main theorem. In \cite{christ1}, Christ has derived an estimate of the Bergman kernel associated to the trivial line bundle defined over  $\C^n$, away from the diagonal. In \cite{delin}, Delin has derived a similar estimate for the  $\mathcal{C}^1$-seminorm of the Bergman kernel associated to the trivial line bundle defined over  $\C^n$, away from the diagonal. 

\vspace{0.1cm}
Let $X$ be a compact K\"ahler manifold, and let $\L$ be a holomorphic line bundle defined over $X$. Then, for any $k\in\N$, off-diagonal estimates of the Bergman kernel associated to $\L^{\otimes k}$, were derived by Christ in \cite{christ2}. Using Szeg{\'{o}} kernels and Poincar{\'e} series, in \cite{lu-zelditch}, Lu and Zelditch have also derived an estimate of the Bergman kernel associated to $\L^{\otimes k}$, away from the diagonal.

\vspace{0.1cm}
Let $X$ be a sympletic manifold of real dimension-$2n$ with bounded geometry, and let $\L$ and $\mathcal{F}$ be a holomorphic line bundle and a vector bundle defined over $X$, respectively. In \cite{ma}, for $n,k\in \N$, Ma and Marinescu have derived estimates of  $\mathcal{C}^{n}$-norms of the Bergman kernel associated to the vector bundle $\mathcal{L}^{\otimes k}\otimes\mathcal{F}$, away from the diagonal. 

\vspace{0.1cm}
When $X$ is a compact hyperbolic Riemann surface, in \cite{chen}, Chen and Fu have derived an estimate of the Bergman kernel associated to the cotangent bundle along the diagonal. In \cite{abms},  the authors have derived an estimate of the Bergman kernel associated to tensor-powers of the line bundle of holomorphic cusp forms of weight-$2$, along the diagonal. As the line bundle of holomorphic cusp forms of weight-$2$ is isometric to the cotangent bundle, the estimate computed in \cite{abms} can also be viewed as an estimate of the Bergman kernel associated to the tensor-powers of the cotangent bundle. The estimates derived in \cite{chen} and \cite{abms} are stable in covers of compact hyperbolic Riemann surfaces. 

\vspace{0.1cm}
When $X$ is a noncompact hyperbolic Riemann surface of finite volume, an estimate of the Bergman kernel associated to the tensor-powers of the line bundle of holomorphic cusp forms of weight-$2$ is obtained in \cite{jk2}, which is also stable in covers. 

\vspace{0.2cm}
{\bf{Statement of main theorem.}}
We now state the main theorem of the article.
\begin{theorem*}
Let $X$ be a hyperbolic Riemann surface of finite volume, and for any $k\in\N$, let $\bkx$ denote the Bergman kernel associated to $\Ok$, where $\O$ denotes the cotangent bundle of $X$. Let $\|\cdot\|_{\mathrm{hyp}}$ denote the Hermitian metric on $\Ok$. For any  $\delta\geq \rx$, let $z=x+iy,\,w=u+iv\in X$ with $\dh(z,w)\geq \delta$, where  $\dh(z,w)$ denotes the geodesic distance between the points $z$ and $w$ on $X$. When $X$ is compact, we have the following estimate
 \begin{align}\label{mainthmeqn1}
\|\bkx\|_{\mathrm{hyp}}(z,w)\leq \cx;
 \end{align}
 when $X$ is noncompact, we have the following estimate 
  \begin{align}\label{mainthmeqn2}
\|\bkx\|_{\mathrm{hyp}}(z,w)\leq  \cx+ \frac{2k-1}{4\pi\cosh^{2k}(\delta\slash 2)}+\frac{(2k-1)\sqrt{yv}}{\sqrt{\pi}}\cdot\frac{\Gamma\big((k-1)\slash 2\big)}{\Gamma(k)},
 \end{align}
where
\begin{align}\label{defncx}
&\cx:=\frac{(2k-1)\sinh(\delta+\rx)}{4\pi\cosh^{2k}\big((\delta-\rx)\slash2\big)\sinh(\rx)}+\frac{(2k-1)\sinh(\delta)}{2\pi\cosh^{2k}(\delta\slash 2)}\cdot\frac{\cosh(\rx\slash 4)}{\sinh(\rx\slash 4)}+\notag\\[0.1cm]&\frac{2k-1}{2\pi(2k-2)\cosh^{2k-2}(\delta\slash2)}\bigg(2+\frac{1}{\sinh^{2}(\rx\slash 4)}\bigg)+\frac{2k-1}{\pi(2k-4)\cosh^{2k-4}(\delta\slash2)}\cdot\frac{1}{\sinh^{2}(\rx\slash 4)}.
\end{align}

\vspace{0.15cm}
 Here $\rx$ denotes the injectivity radius of $X$, which is as defined in equations \eqref{irad1} and \eqref{irad2}, for compact and noncompact hyperbolic Riemann surfaces, respectively.  
\end{theorem*}

\vspace{0.15cm}
\begin{rem}
Similar to the estimates of the Bergman kernel derived in \cite{chen}, \cite{abms}, and \cite{jk2}, it is easy to show that our estimates \eqref{mainthmeqn1} and \eqref{mainthmeqn2} are stable in covers of Riemann surfaces, by following similar arguments as in \cite{abms}. In the noncompact setting, our off-diagonal estimate of the Bergman kernel  \eqref{mainthmeqn2} is the same as the estimate of the Bergman kernel along the diagonal, which is derived in \cite{jk2}. 
\end{rem}

{\bf{Acknowledgements.}} Both the authors acknowledge the support of INSPIRE research grant DST/INSPIRE/04/2015/002263.
\section{Background material}
In this section, we set up the notation, and recall the background details needed for the proofs of the Main theorem. 

\vspace{0.15cm}
 Let $X$ denote a hyperbolic Riemann surface of finite volume. By uniformization theorem from complex analysis, $X$ can be realized as the quotient space 
$\G\backslash\H$, where $\G\subset\mathrm{PSL}_{2}(\R)$ is a cofinite Fuchsian subgroup, and $\H$ is the complex upper half-plane. Locally, we identify $X$ with its universal cover $\H$. 

\vspace{0.2cm}
When $X$ is a noncompact hyperbolic Riemann surface of finite volume, without loss of  generality, we assume that the point $i\infty$ is the only puncture of $X$, which we denote by $\infty$. 

\vspace{0.15cm}
Let $\hyp$ denote the natural hyperbolic metric on $X$, which is of constant negative curvature $-1$. For any $z,w\in X$, let $\dh(z,w)$ denote the geodesic distance on the Riemann surface $X$, which is the natural distance function on $X$, coming from the hyperbolic metric $\hyp$. 

\vspace{0.15cm}
The injectivity radius of a compact hyperbolic Riemann surface $X$ is given by the following formula
\begin{align}\label{irad1}
\rx:=\inf\big\lbrace \dh(z,\gamma z)|\,z\in X,\, \gamma\in\Gamma\backslash\lbrace\mathrm{Id}\rbrace\big\rbrace;
\end{align}
when $X$ is a noncompact hyperbolic Riemann surface of finite volume, we define the injectivity radius of $X$ by the following formula
\begin{align}\label{irad2}
\rx:=\inf\big\lbrace \dh(z,\gamma z)|\,z\in X,\, \gamma\in\Gamma\backslash\Gamma_{\infty}\rbrace,
\end{align}
where $\Gamma_{\infty}$ is the stabilizer of the cusp $\infty$. 

\vspace{0.15cm}
Let $\O$ denote the cotangent bundle of holomorphic differential $1$-forms on $X$. The global sections of this line bundle are of the form $f(z)dz$, where $f(z)$ is a holomorphic cusp form of weight-$2$ with respect to $\Gamma$. So $H^{0}(X,\O)$, the space of global sections of $\O$ can be identified with $\Stwo$, the complex vector-space of weight-$2$ cusp forms.  This implies that for any $k\in\N$, we can identify $H^{0}(X,\Ok)$ with $\Sk$ as complex vector-spaces.

\vspace{0.15cm}
For any $k\in\N$, let $\omega\in H^{0}(X,\Omega_{X}^{\otimes k})$ be a global differential $k$-form. Then, locally, at any $z\in X$, $\omega(z)=f(z)dz^{\otimes k}$, where $f$ is a weight-$2k$ cusp form with respect to $\Gamma$. Furthermore, there exists a point-wise metric on $H^{0}(X,\Ok)$, which is denoted by $\|\cdot\|_{\mathrm{hyp}}$, and  at the point $z=x+iy\in X$, it is given by the following formula
\begin{align*}
\|\omega\|_{\mathrm{hyp}}(z)=y^{k}|f(z)|.
\end{align*}

\vspace{0.15cm}
For any $\omega(z)=f(z)dz^{\otimes k},\eta(z)=g(z)dz^{\otimes k}\in H^{0}(X,\Ok)$, the $L^{2}$-metric induced by  $\|\cdot\|_{\mathrm{hyp}}$ on $H^{0}(X,\Ok)$ is denoted by $\langle\cdot,\cdot\rangle_{\mathrm{hyp}}$, and is given by the following formula
\begin{align*}
 \langle \omega,\eta\rangle_{\mathrm{hyp}}=\int_{X}y^{2k}f(z)\overline{g}(z)\hyp(z).
\end{align*}

Let $\omega_{1},\ldots,\omega_{\jk}$ denote an orthonormal basis of $H^{0}(X,\Ok)$ with respect to the $L^2$-metric $\langle\cdot,\cdot\rangle_{\mathrm{hyp}}$. Then, for any $z,w\in X$,  the Bergman kernel $\bkx(z,w)$ associated to the line bundle $\Ok$ is given by the following formula
\begin{align*}
\bkx(z,w)=\sum_{i=1}^{\jk}\omega_{i}(z)\overline{\omega}_{i}(w).
\end{align*}
It is easy to show that the Bergman kernel is independent of the choice of the orthonormal basis of $H^{0}(X,\Ok)$. 

\vspace{0.2cm}
We now describe the Bergman kernel associated to $\Sk$, the vector space of weight-$2k$ cusp forms. For $f\in\Sk$, we have the following point-wise metric at $z=x+iy\in X$ 
\begin{align*}
\|f\|_{\mathrm{pet}}(z):=y^{k}|f(z)|,
\end{align*}
which induces an $L^{2}$-metric on $\Sk$, which is also known as the Petersson inner-product. For any $f,g\in\Sk$, the Petersson inner-product is given by the following formula
\begin{align*}
\langle f,g\rangle_{\mathrm{pet}}:=\int_{X}y^{2k}f(z)\overline{g}(z)\hyp(z).
\end{align*}
Let $\lbrace f_1,\ldots,f_{j_{k}}\rbrace $ denote an orthonormal basis for $\Sk$ with respect to the  Petersson inner-product. Then, for any $z=x+iy, w=u+iv\in X$, the Bergman kernel associated to the complex vector-space $\Sk$ is given by the following formula
\begin{align*}
\bk(z,w):=\sum_{i=1}^{j_{k}}f_{i}(z)\overline{f}_{i}(w).
\end{align*}
The Bergman kernel $\bk(z,w)$ is a holomorphic cusp form of weight-$2k$ in $z$, and an anti-holomorphic cusp form of weight-$2k$ in $w$. It can also be given by the following infinite series (see Proposition 1.3 on p. 77 in \cite{frietag})
\begin{align*}
&\bk(z,w)= \frac{(2k-1)(2i)^{2k}}{4\pi}\sum_{\gamma\in\Gamma}\frac{1}{(\gamma z-\overline{w})^{2k}}\cdot\frac{1}{j(\gamma,z)^{2k}},\\
&\mathrm{where}\,\,\mathrm{for\,\,any}\,\gamma=\left(\begin{array}{cc} a&b\\c&d\end{array}\right)\in\Gamma, \,j(\gamma,z)=cz+d.
\end{align*}
As $H^{0}(X,\Ok)\cong\Sk$ as complex vector-spaces, we have the following relation of Bergman kernels
\begin{align*}
\bkx(z,w)=\bk(z,w)\big(dz^{\otimes k}\wedge d\overline{w}^{\otimes k}\big),
\end{align*}
and the point-wise metric on $\Ok$ induces the following point-wise metric on $\bkx(z,w)$
\begin{align}\label{seriesbergmanreln0}
\|\bkx\|_{\mathrm{hyp}}(z,w)=y^{k}v^{k} \frac{(2k-1)(2)^{2k}}{4\pi}\Bigg|\sum_{\gamma\in\Gamma}\frac{1}{(\gamma z-\overline{w})^{2k}}\cdot\frac{1}{j(\gamma,z)^{2k}}\Bigg|.
\end{align}

\vspace{0.15cm}
Recall that for any $z=x+iy,\,w=u+iv\in X$, we have the following formula
\begin{align*}
\cosh^{2}\big(\dh(z,w)\slash 2\big)=\frac{|z-\overline{w}|^{2}}{4yv}.
\end{align*}
Using the above relation and equation \eqref{seriesbergmanreln0}, we derive the following inequality
\begin{align}\label{seriesbergmanreln}
\|\bkx\|_{\mathrm{hyp}}(z,w)\leq \frac{2k-1}{4\pi}\sum_{\gamma\in\Gamma}\Bigg|\frac{4^ky^{k}v^{k} }{(\gamma z-\overline{w})^{2k}j(\gamma,z)^{2k}} \Bigg|= \notag\\\frac{2k-1}{4\pi}\sum_{\gamma\in\Gamma}\frac{\big(4\mathrm{Im}(\gamma z)v\big)^{k} }{\big|\gamma z-\overline{w}\big|^{2k}} = \frac{2k-1}{4\pi}\sum_{\gamma\in\Gamma}\frac{1}{\cosh^{2k}\big(\dh(\gamma z,w)\slash 2\big)}.
\end{align}

\vspace{0.15cm}
For a hyperbolic Riemann surface of finite volume, we now state an inequality from \cite{jl}, which is adapted to our setting. The inequality gives us an estimate for the number of elements in $\G$ or $\G\backslash \G_{\infty}$, depending on whether $X$ is compact or noncompact, respectively.

\vspace{0.15cm}
 For any positive, smooth, real-valued, and decreasing function $f$ defined on $\mathbb{R}_{\geq 0}$, and for any $\delta > \rx\slash 2$, we have the following inequality 
\begin{align}\label{jlineq1}
 \int_0^{\infty} f(\rho) dN_{\G}(z,w;\rho) \leq \int_0^{\delta} f(\rho) dN_{\Gamma}(z,w;\rho) + f(\delta) \frac{2\cosh(\rx\slash 4)\sinh(\delta)}{\sinh(\rx\slash 4)}
 +\notag\\\frac{1}{2\sinh^2(\rx\slash 4)} \int_{\delta}^{\infty} f(\rho) \sinh(\rho +\rx\slash 2) d\rho,
\end{align}
where 
\begin{align*}
N_{\Gamma}(z,w;\rho) := \mathrm{card}\,\{ \gamma |\,\gamma \in \G\backslash\Gamma_{\infty},\,\dh(\gamma z,w) \leq \rho \}.
\end{align*}  

\vspace{0.15cm}
From arguments similar to the ones used in deriving inequality \eqref{jlineq1} in \cite{jl}, for any $\delta>0$, and $z,w\in X$, we have the following inequality
\begin{align}\label{jlineq2}
N_{\G}(z,w;\delta)\leq \frac{\sinh(\delta+\rx)}{\sinh(\rx)}.
\end{align}
The above inequality has already been used in the above form in \cite{jkcomp}. Notice that our definition for injectivity radius is two times the injectivity radius in \cite{jl}, and both our inequalities, \eqref{jlineq1} and \eqref{jlineq2} take this fact into account. 

\vspace{0.15cm}
Here, it is understood that, when $X$ is compact, $\G_{\infty}=\emptyset$. So inequalities \eqref{jlineq1} and \eqref{jlineq2} also hold true in the compact setting. 
\section{Proof of the Main theorem}
\begin{proof}[Proof of estimate \eqref{mainthmeqn1}]
Given that $X=\G\backslash\H$ is a compact hyperbolic Riemann surface with injectivity radius $\rx$, and $\delta\geq \rx$. Furthermore, $\bkx(z,w)$ is the Bergman kernel for the line bundle $\Ok$.  Combining inequalities \eqref{seriesbergmanreln}  and \eqref{jlineq1}, for any $z,w\in X$ with $\dh(z,w)\geq \delta$, we find
\begin{align}\label{proof1eqn1}
&\|\bkx\|_{\mathrm{hyp}}(z,w)\leq \frac{2k-1}{4\pi}\sum_{\gamma\in\Gamma}\frac{1}{\cosh^{2k}\big(\dh(\gamma z,w)\slash 2\big)}=\notag\\[0.1cm]&\frac{2k-1}{4\pi}\int_{0}^{\delta}\frac{dN_{\G}(z,w;\rho)}{\cosh^{2k}\big(\dh(\gamma z,w)\slash 2\big)}+\frac{(2k-1)\sinh(\delta)}{2\pi\cosh^{2k}(\delta\slash2)}\cdot
\frac{\cosh(\rx\slash 4)}{\sinh(\rx\slash 4)} +\notag\\[0.1cm]&\frac{2k-1}{8\pi\sinh^2(\rx\slash 4)}\int_{\delta}^{\infty}\frac{\sinh(\rho +\rx\slash 2) d\rho}{\cosh^{2k}(\rho\slash2)}.
\end{align}

\vspace{0.1cm}
We now estimate the first term on the right-hand side of the equality in the above inequality. For  $\gamma\in \G$, and $z,w\in X$ with $\dh(z,w)\geq \delta$, using triangular inequality, we derive  
\begin{align*}
\dh(z,\gamma z)+\dh(\gamma z, w)\geq \dh(z,w)\geq \delta.
\end{align*} 
Using which, we compute 
\begin{align*}
\inf_{\gamma\in\G\backslash\mathrm{Id}}\big(\dh(z,\gamma z)+\dh(\gamma z, w)\big)&=\inf_{\gamma\in\G\backslash\mathrm{Id}}\dh(z,\gamma z)+\inf_{\gamma\in\G\backslash\mathrm{Id}}\dh(\gamma z,w)\\&=\rx+\inf_{\gamma\in\G\backslash\mathrm{Id}}\dh(\gamma z, w)\geq \delta,
\end{align*} 
which implies that for any $\gamma\in \G$, we have
\begin{align*}
\dh(\gamma z, w)\geq \delta-\rx\geq 0\implies  \frac{1}{\cosh^{2k}\big((\delta-\rx)\slash 2\big)}\geq 
\frac{1}{\cosh^{2k}\big(\dh(\gamma z,w)\slash 2\big)},
\end{align*}
for any $k\in\N$. So combining the above inequality with inequality \eqref{jlineq2}, we arrive at the following inequality
\begin{align*}
&\int_{0}^{\delta}\frac{dN_{\G}(z,w;\rho)}{\cosh^{2k}\big(\dh(\gamma z,w)\slash 2\big)}\leq \notag\\ &\big(N_{\G}(z,w;\delta)\big)\cdot\sup_{\gamma\in N_{\G}(z,w;\delta)}\bigg(\frac{1}{\cosh^{2k}\big(\dh(\gamma z, w)\slash 2\big)}\bigg)\leq\frac{\sinh(\delta+\rx)}{\cosh^{2k}\big((\delta-\rx)\slash2\big)\sinh(\rx)},
\end{align*}
which implies that we have the following estimate for the first term on the right-hand side of the equality in \eqref{proof1eqn1}
\begin{align}\label{proof1eqn2}
\frac{2k-1}{4\pi}\int_{0}^{\delta}\frac{dN_{\G}(z,w;\rho)}{\cosh^{2k}\big(\dh(\gamma z,w)\slash 2\big)}\leq\frac{(2k-1)\sinh(\delta+\rx)}{4\pi\cosh^{2k}\big((\delta-\rx)\slash2\big)\sinh(\rx)}.
\end{align}

\vspace{0.1cm}
We now estimate the third term on the right-hand side of the equality in \eqref{proof1eqn1}. For any $\rho\geq \delta$, observe that
\begin{align*}
&\sinh(\rho+\rx\slash 2)=\sinh(\rho)\cosh(\rx\slash 2)+\cosh(\rho)\sinh(\rx\slash 2)\leq\\[0.1cm]& \sinh(\rho)\cosh(\rx\slash 2)+\cosh(\rho)\sinh(\rho)=
2\sinh(\rho\slash 2)\cosh(\rho\slash 2)\big(\cosh(\rx\slash 2)+\cosh(\rho)\big)\leq \\[0.1cm]&2\sinh(\rho\slash 2)\cosh(\rho\slash 2)\big(\cosh(\rx\slash 2)+2\cosh^{2}(\rho\slash 2)\big).
\end{align*}

\vspace{0.1cm}
Using which, we derive that
\begin{align*}
&\int_{\delta}^{\infty}\frac{\sinh(\rho +\rx\slash 2) d\rho}{\cosh^{2k}(\rho\slash2)}\leq \cosh(\rx\slash 2)\int_{\delta}^{\infty}\frac{2\sinh(\rho\slash 2)d\rho}{\cosh^{2k-1}(\rho\slash 2)}+\int_{\delta}^{\infty}\frac{4\sinh(\rho\slash 2)d\rho}{\cosh^{2k-3}(\rho\slash 2)}=\notag\\[0.1cm]&\frac{4\cosh(\rx\slash 2)}{(2k-2)\cosh^{2k-2}(\delta\slash2)}+\frac{8}{(2k-4)\cosh^{2k-4}(\delta\slash2)}.
\end{align*}

\vspace{0.1cm}
So we have the following estimate for the third term of right-hand side of the equality in \eqref{proof1eqn1}
\begin{align}\label{proof1eqn3}
&\frac{2k-1}{8\pi\sinh^2(\rx\slash 4)}\int_{\delta}^{\infty}\frac{\sinh(\rho +\rx\slash 2) d\rho}{\cosh^{2k}(\rho\slash2)}
\leq\notag \\[0.1cm]&\frac{(2k-1)}{2\pi(2k-2)\cosh^{2k-2}(\delta\slash2)} \cdot\frac{\cosh(\rx\slash 2)}{ \sinh^2(\rx\slash 4)} +\frac{2k-1}{\pi (2k-4)\cosh^{2k-4}(\delta\slash2)}\cdot\frac{1}{\sinh^2(\rx\slash 4)}=\notag\\[0.1cm] &\frac{(2k-1)}{2\pi(2k-2)\cosh^{2k-2}(\delta\slash2)} \bigg(2+\frac{1}{\sinh^{2}(\rx\slash 4)}\bigg)+\frac{2k-1}{\pi (2k-4)\cosh^{2k-4}(\delta\slash2)}\cdot\frac{1}{\sinh^2(\rx\slash 4)}.
\end{align}

\vspace{0.1cm}
The proof of estimate \eqref{mainthmeqn1} follows from combining estimates \eqref{proof1eqn1}, \eqref{proof1eqn2}, and \eqref{proof1eqn3}.
\end{proof}

\vspace{0.2cm}
\begin{proof}[Proof of estimate \eqref{mainthmeqn2}]
Now let $X$ be a noncompact hyperbolic Riemann surface of finite volume. For $z,w\in X$ with $\dh(z,w)\geq \delta$, from inequality \eqref{seriesbergmanreln}, we have
\begin{align}\label{proof2eqn1}
&\|\bkx\|_{\mathrm{hyp}}(z,w)\leq \frac{2k-1}{4\pi}\sum_{\gamma\in\G\backslash\Gamma_{\infty}}\frac{1}{\cosh^{2k}(\dh(z,w)\slash 2)}+\frac{2k-1}{4\pi}\sum_{\gamma\in\Gamma_{\infty}}\frac{1}{\cosh^{2k}\big(\dh(z,w)\slash 2\big)}.
\end{align}
From similar arguments as in the proof of estimate \eqref{mainthmeqn1}, we have the following estimate for the first term on the right-hand side of the above inequality
\begin{align}\label{proof2eqn2}
\frac{2k-1}{4\pi}\sum_{\gamma\in\G\backslash\Gamma_{\infty}}\frac{1}{\cosh^{2k}(\dh(z,w)\slash 2)} \leq \cx,
\end{align}
where $\cx$ is as in equation \eqref{defncx}.

\vspace{0.1cm}
We now estimate the second term on the right-hand side of inequality \eqref{proof2eqn1}. Without loss of generality, we assume that 
\begin{align*}
\Gamma_{\infty}:=\bigg\lbrace\left(\begin{array}{cc} 1&n\\0&1\end{array}\right)\big|\,n\in\mathbb{Z} \bigg\rbrace.
\end{align*}

\vspace{0.1cm}
This implies that for $z=x+iy, w=u+iv\in X$, we have
\begin{align}\label{proof2eqn3}
&\sum_{\gamma\in\Gamma_{\infty}}\frac{1}{\cosh^{2k}\big(\dh(\gamma z,w)\slash 2\big)}=\frac{1}{\cosh^{2k}\big(\dh(z,w)\slash 2\big)}+\sum_{n\in\mathbb{Z}\backslash\lbrace0\rbrace}\frac{(4yv)^{k}}{\big((x+n-u)^{2}+(y+v)^{2}\big)^{k}}\leq \notag\\[0.1cm]&\frac{1}{\cosh^{k}(\delta\slash 2)}+\int_{0}^{\infty}\frac{d\alpha}{\bigg(\big(\frac{\alpha+x-u}{2\sqrt{yv}}\big)^{2}+\frac{(y+v)^{2}}{4yv}\bigg)^{k}}+\int_{0}^{\infty}\frac{d\alpha}{\bigg(\big(\frac{\alpha+u-x}{2\sqrt{yv}}\big)^{2}+\frac{(y+v)^{2}}{4yv}\bigg)^{k}}.
\end{align}

\vspace{0.1cm}
Using the fact that for $y,v >0$, $(y+v)^{2}\geq 4yv$, for the first integral on the right-hand side of the above inequality, we compute
\begin{align}\label{proof2eqn4}
\int_{0}^{\infty}\frac{d\alpha}{\bigg(\big(\frac{\alpha+x-u}{2\sqrt{yv}}\big)^{2}+\frac{(y+v)^{2}}{4yv}\bigg)^{k}}\leq \int_{0}^{\infty}\frac{d\alpha}{\bigg(\big(\frac{\alpha+x-u}{2\sqrt{yv}}\big)^{2}+1\bigg)^{k}}.
\end{align}

\vspace{0.1cm}
Making the substituting $(\alpha+x-u)\slash2\sqrt{yv}=\theta$, and using formula 3.251.2 from \cite{gr}, we arrive at the following inequality
\begin{align}\label{proof2eqn5}
 \int_{0}^{\infty}\frac{d\alpha}{\bigg(\big(\frac{\alpha+x-u}{2\sqrt{yv}}\big)^{2}+1\bigg)^{k}}\leq \int_{-\infty}^{\infty}\frac{2\sqrt{yv}d\theta}{\big(\theta^{2}+1\big)^{k}}=\frac{4\sqrt{\pi yv}\G\big(k-1\slash 2\big)}{\G(k)}.
\end{align}
Following similar arguments, we have the following estimate for the second integral on the right-hand side of inequality \eqref{proof2eqn3}
\begin{align}\label{proof2eqn6}
\int_{0}^{\infty}\frac{d\alpha}{\bigg(\big(\frac{\alpha+u-x}{2\sqrt{yv}}\big)^{2}+\frac{(y+v)^{2}}{4yv}\bigg)^{k}}\leq \frac{4\sqrt{\pi yv}\G\big(k-1\slash 2\big)}{\G(k)}.
\end{align}

\vspace{0.1cm}
Combining inequalities \eqref{proof2eqn3}, \eqref{proof2eqn5}, and \eqref{proof2eqn6}, we arrive at the following estimate for the second term on the right-hand side of 
inequality \eqref{proof2eqn1}
\begin{align}\label{proof2eqn7}
\frac{2k-1}{4\pi}\sum_{\gamma\in\Gamma_{\infty}}\frac{1}{\cosh^{2k}\big(\dh(z,w)\slash 2\big)}\leq\frac{2k-1}{4\pi\cosh^{2k}(\delta\slash 2)}+\frac{(2k-1)\sqrt{yv}}{\sqrt{\pi}}\cdot\frac{\Gamma\big((k-1)\slash 2\big)}{\Gamma(k)}.
\end{align}
Combining estimates \eqref{proof2eqn1}, \eqref{proof2eqn2}, and \eqref{proof2eqn7} completes the proof of estimate \eqref{mainthmeqn2}, and also the proof of Main theroem.
\end{proof}
\begin{rem}
Let $X$ be  hyperbolic Riemann surface of finite volume with injectivity radius $\rx$. When $X$ is compact, for any $\delta\geq \rx$, and $z=x+iy,\,w=u+iv\in X$ with $\dh(z,w)\geq \delta$, a careful analysis of each of the term comprising the constant $\cx$ given in equation \eqref{defncx}, leads us to the conclusion that 
\begin{align*}
\cx=O_{X}\bigg(\frac{k}{\cosh^{2k-4}\big((\delta-\rx)\slash 2\big)}\bigg).
\end{align*}
Using the fact that $\cosh(u)\geq e^u\slash 2$, for all $u\geq 0$, we observe that
\begin{align*}
\cx=O_{X}\Bigg(\frac{k\cdot2^{2k-4}}{e^{(k-2)(\delta-\rx)\big)}}\Bigg)=O_{X}\bigg(ke^{-(k-2)\big(\delta-\rx-2\ln2\big)}\bigg).
\end{align*}
For $\delta\geq \rx$ and $k\in \N$ sufficiently large, our estimate \eqref{mainthmeqn1} is a slightly stronger estimate than the more general estimates derived in \cite{ma} and \cite{christ2}.

\vspace{0.1cm}
When $X$ is noncompact, adapting the same arguments as in p. 12 in section 5 of \cite{jk2}, we find that
\begin{align*}
\frac{(2k-1)\sqrt{yv}}{\sqrt{\pi}}\cdot\frac{\Gamma\big((k-1)\slash 2\big)}{\Gamma(k)}=O\big(k^{3\slash 2}\big),
\end{align*}
which implies that
\begin{align*}
\|\bkx\|_{\mathrm{hyp}}(z,w)=O_{X}\big(k^{3\slash 2}\big),
\end{align*}
which is the same as the estimate for the Bergman kernel along the diagonal derived in \cite{jk2}. 
\end{rem}

\vspace{0.15cm}
\begin{rem}
Let $X_1=\G_1\backslash\H$, $X_{0}=\G_{0}\backslash\H$ be two compact hyperbolic Riemann surfaces. Let $X_1$ be a finite cover of $X_0$, which implies that $\G_1$ is a finite index subgroup of $\G_{0}$. So for any $k\in\N$, and $z,w \in X_{1}$ with $\dh(z,w)\geq \delta$, from estimates \eqref{seriesbergmanreln} and \eqref{mainthmeqn1} we find that
\begin{align*}
&\|\mathcal{B}^{k}_{\Omega_{X_1}}\|_{\mathrm{hyp}}(z,w)\leq  \frac{2k-1}{4\pi}\sum_{\gamma\in\Gamma_1}\frac{1}{\cosh^{2k}\big(\dh(\gamma z,w)\slash 2\big)}\leq  \\&\frac{2k-1}{4\pi}\sum_{\gamma\in\Gamma_0}\frac{1}{\cosh^{2k}\big(\dh(\gamma z,w)\slash 2\big)}\leq \mathcal{C}_{X_0}.
\end{align*}

\vspace{0.1cm}
This implies that our estimate \eqref{mainthmeqn1} is stable in covers of compact hyperbolic Riemann  surfaces. Following the same argument, we can conclude that our estimate \eqref{mainthmeqn2} is also stable in covers of noncompact hyperbolic Riemann surfaces of finite volume. 
\end{rem}
{}

\vspace{0.3cm}
{\small{Anilatmaja Aryasomayajula,\\
Department of Mathematics,\\
Indian Institute of Science Education and Research Tirupati, \\
Karkambadi Road, Mangalam (B.O),Tirupati-517507, India.\\
email: {\it{anil.arya@iisertirupati.ac.in}}

\vspace{0.2cm}
Priyanka Majumder, \\
Department of Mathematics, \\
Indian Institute of Science Education and Research Tirupati, \\
Karkambadi Road, Mangalam (B.O),Tirupati-517507, India.\\
email: {\it{pmpriyanka57@gmail.com }}}}
\end{document}